\documentclass[12pt,reqno,oneside]{amsart}
 \usepackage{verbatim,amsmath,amssymb,cite,epsfig,xspace}
\usepackage{color,cite,graphicx,esint}

\numberwithin{equation}{section}

\theoremstyle{plain}
\newtheorem{theorem}{Theorem}[section]
\newtheorem{lemma}{Lemma}[section]
\newtheorem{corollary}{Corollary}[section]
\theoremstyle{definition}
\newtheorem{definition}{Definition}[section]

\theoremstyle{remark}

\begin{document}
\title[]{Local-in-space estimates near initial time for weak solutions of the Navier-Stokes equations and forward self-similar solutions}
\author{Hao Jia, Vladim\'{i}r \v{S}ver\'{a}k \\ \\ {\smaller University of Minnesota}}
\maketitle

\begin{center}
\textbf{Abstract}\end{center}
We show that the classical Cauchy problem for the incompressible 3d Navier-Stokes equations with $(-1)$-homogeneous initial data has a global scale-invariant solution which is smooth for positive times. Our main technical tools are local-in-space regularity estimates near the initial time, which are of independent interest.

\begin{section}{Introduction}
We consider the classical Cauchy problem for the incompressible Navier-Stokes equation in $R^3\times (0,\infty)$
\begin{eqnarray}
\left.\begin{array}{rcl}
u_t+u\nabla u +\nabla p -\Delta u & = & 0 \\
{\rm div \,}u & = & 0
\end{array} \right\}& &  \quad \hbox{in $R^3\times(0,\infty)$}\,, \label{nse1} \\
u|_{t=0}\,\,\, = \,\,\,u_0\,\,\,\, && \quad \hbox{in $R^3$}\,.\label{nse2}
\end{eqnarray}
We recall that the problem is invariant under the scaling
\begin{equation}\label{scaling}
\begin{array}{rcccl}
u(x,t) & \to & u_\lambda(x,t) & = & \lambda u(\lambda x, \lambda ^2t)\,,\\
p(x,t) & \to & p_\lambda(x,t) & = & \lambda^2p(\lambda x,\lambda^2 t)\,,\\
u_0(x) &\to & u_{0\lambda} (x) & = & \lambda u_0(\lambda x)\,,
\end{array}
\end{equation}
where $\lambda>0$. We say that a solution $u$ is {\it scale-invariant} if $u_\lambda=u$ and $p_\lambda=p$ for each $\lambda>0$. Similarly, we say that an initial condition $u_0$ is scale-invariant, if $u_{0\lambda}=u_0$ for each $\lambda>0$. This is of course the same as requiring that $u_0$ be $(-1)-$ homogeneous.

One of our goals in this paper is to give a proof of the following result.

\begin{theorem}\label{thm1}
Assume $u_0$ is scale-invariant and locally H\"older continuous in $R^3\setminus\{0\}$ with ${\rm ~div~~}u_0=0$ in $R^3$. Then the Cauchy problem~(\ref{nse1}),~(\ref{nse2}) has at least one scale-invariant solution $u$ which is smooth in $R^3\times(0,\infty)$ and locally H\"older continuous in $R^3\times[0,\infty)\setminus \{(0,0)\}$.
\end{theorem}
Previously this result has been known only under suitable smallness conditions on $u_0$, see for example~\cite{CaPl, KoTa}. For small $u_0$ one can also prove uniqueness (in suitable function classes). It is quite conceivable that uniqueness may fail for large data. We will comment on this point in more detail below.

\smallskip
The second important theme of our paper can be perhaps called local-in-space regularity estimates near the initial time $t=0$. It is known that if $u_0\in L^q$ for $q\ge 3$, then the initial value problem~(\ref{nse1}),~(\ref{nse2}) has a unique local-in-time ``mild solution" defined on some time interval $(0,T)$, which is smooth in $R^3\times(0,T)$ and has in many respects the same regularity as the solution of the heat equation in $R^3\times[0,T)$ for times close to $t=0$, see for example~\cite{DD, GPS}. A natural question is under which condition this result can be localized in space: if $u_0$ is a quite general initial condition  for which a generalized suitable weak solution $u$ in the sense of~\cite{LemRie} is defined and $u_0|_{B_r}\in L^q(B_r)$ for some $q>3$, say, can we conclude that $u$ is regular in $B_{r\over 2}\times[0,t_1)$ for some time $t_1>0$?
We prove that this is indeed the case under quite general assumptions, which include $u_0$ which is in $L^2_{\rm loc}$ and $\int_{B_{x,r}} |u_0|^2\,dx\to 0$ for $x\to\infty$. Due to non-local effects of the pressure the solution $u$ in $B_{r\over 2}\times[0,t_1)$ may not have the same amount of regularity as the solution of the heat equation in this situation, but the non-local effects are limited to the influence of the ``harmonic part of the pressure" in a suitable pressure decomposition. We can formulate this type of results somewhat loosely in the following statement.

\bigskip
\noindent
{\bf (S)$\,\,$}{\sl Modulo the usual (and quite mild) non-local influences of the pressure, local regularity of the initial data propagates for at least a short time.}

\bigskip
We refer the reader to Section 3 for precise statements. Statement (S), in addition to being of independent interest, is one of the main ingredients of our proof of Theorem~\ref{thm1}.

\smallskip
Results in the direction of (S) can be found already in the classical paper~\cite{CKN}. More recently, related questions about vorticity propagation have been studied in~\cite{Tao}. Our main result concerning (S), Theorem~\ref{th:localRegularityForLeraySolution}, takes a somewhat different angle on (S).

Our proof of (S) (see also Theorem~\ref{th:localRegularityForLeraySolution}) is based on a combination of techniques from~\cite{HaSv, Lin, LaSe, LemRie, Seregin}. Heuristically, the main point is that one can obtain a sufficient control of the energy flux into ``good regions" from the rest of the space, see Section 3. Once we know that only small amount of energy can move into the ``good region" one can use (a slight modification of) partial regularity schemes in~\cite{Lin,LaSe} to prove regularity.

\smallskip
To prove Theorem 1.1, we seek the solution $u(x,t)$ in the form
\begin{equation}\label{ansatz}
u(x,t)={1\over \sqrt{t}}\,U\left({x\over \sqrt t}\right)\,.
\end{equation}
The Navier-Stokes equation for $u$ gives
\begin{equation}\label{Ueq}
-\Delta U -{1\over 2}U-{1\over 2}x\nabla U +U\nabla U+\nabla P=0,\qquad {\rm div\,} U=0\,,
\end{equation}
in $R^3$. For a scale-invariant $u_0$ the problem of finding a scale-invariant solution of the Cauchy problem~(\ref{nse1}),~(\ref{nse2}) is equivalent to the problem of finding a solution of~(\ref{Ueq}) with the asymptotics
\begin{equation}\label{asymptotics}
|U(x)-u_0(x)|=o\left({1\over |x|}\right)\,,\qquad x\to\infty\,.
\end{equation}
The problem~(\ref{Ueq}),~(\ref{asymptotics}) is reminiscent of the classical Leray's problem of finding steady-state solution of the Navier-Stokes equation in a bounded domain, with given boundary conditions for $U$. We will show that one can solve this problem using the Leray-Schauder degree theory, just as in the case of the bounded domain. The non-trivial part is to find the right functional-analytic setup and establish the necessary a-priori estimates. The main difficulty is to find good estimates near $\infty$. This difficulty will be overcome by applying statement (S) above. Heuristically it is clear that when $u$ is given by~(\ref{ansatz}), then estimates of $u$ near $t=0$ are closely related to estimates of $U$ near $\infty$. In Section 4 we will make this more precise.

\smallskip
As in the case of the bounded domains, the Leray-Schauder approach gives existence of the solutions, but not uniqueness. In the case of bounded domains one does not generically expect uniqueness for large data, and this non-uniqueness is in fact expected to be quite typical in the context of the steady Navier-Stokes, once the data is  large. Could this also be the case for the problem~(\ref{Ueq}),~(\ref{asymptotics})? This would lead to non-uniqueness for the Cauchy problem~(\ref{nse1}),~(\ref{nse2}) with scale-invariant $u_0$, and by a suitable truncation of $u_0$ at large $|x|$ possibly also to non-uniqueness for the Leray-Hopf solutions of the Cauchy problem for $u_0\in L^2$. We plan to address these issues in future work.

\medskip

Our paper is organized as follows: in section 2, we prove an `$\epsilon$-regularity' criteria for a generalized Navier-Stokes system; in section 3, we study the local in space near initial time smoothness of Leray solutions; in section 4 we study the asymptotics of forward self similar solutions to Stokes and Navier-Stokes equations; in section 5 we prove the existence of forward self similar solutions for large $-1$ homogeneous initial data.

\smallskip
\noindent
\textbf{Notation:} We use standard notations in our paper. For instance, $B_R(x_0)$ denotes a ball centered at $x_0$ with radius $R$ in $R^3$, $B_R:=B_R(0)$; for $z_0=(x_0,t_0)$, $Q(R,z_0):=B_R(x_0)\times (t_0-R^2,t_0]$, $Q_R:=Q(R,(0,0))$; for any $f$ in $\mathcal{O}$, $\fint_{\mathcal{O}}f:=\frac{1}{|\mathcal{O}|}\int_{\mathcal{O}}f$. We also use the following standard notations in the literature:  for vectors $a$ and $v$, $a\otimes v$ is a matrix with $(a\otimes v)_{ij}=a_iv_j$; for two matices $a,~b$, $(a:b)=a_{ij}b_{ij}$ where we assume the usual Einstein summation convention; for a matrix valued function $f=(f_{ij})$, ${\rm div~~}f$ is a vector with $({\rm div~~}f)_i=(\sum_{j}\partial_jf_{ij})$; $(u)_{Q(R,z_0)}:=\fint_{Q(R,z_0)}u dz$, $(u)_{r}:=(u)_{Q_r}$; $(p)_{B_R(x_0)}(t):=\fint_{B_{R(x_0)}}p(x,t)dx$, $(p)_{r}(t):=(p)_{B_r}(t)$; $Y(u,p,Q(R,z_0)):=(\fint_{Q(R,z_0)}|u-(u)_{Q(R,z_0)}|^3dz)^{1/3}+R(\fint_{Q(R,z_0)}|p-(p)_{B_R(x_0)}(t)|^{3/2}dz)^{2/3}$; $Y(u,p,Q_R):=Y(u,p,Q(R,(0,0)))$. We use $C$ to denote an absolute and often large positive number, $c$ a positive small absolute number, $\epsilon$ the positive small numbers, $C(\alpha,\beta,\dots)$ when the number depends on the parameters $\alpha,~\beta,~\dots$. $C_{{\rm par}}^{\alpha}(\mathcal{O})$ denotes the H\"{o}lder space with respect to the parabolic distance when $\mathcal{O}$ is a space time domain. We adopt the convention that nonessential constants can change from line to line. We use $u_0$ as a divergence free initial data throughout the paper, unless defined otherwise.\\
\begin{section}{ $\epsilon$-regularity criteria}
Our goal in this section is to prove an $\epsilon$-regularity criteria similar to that of Caffarelli-Kohn-Nirenberg for a generalized Navier-Stokes equation. Our setting is as follows:\\
Let $\mathcal{O}$ be an open subset of $R^3_x\times R_t$, $a\in L^m_{loc}(\mathcal{O})$ with $m>5$ (not necessarily an integer), $\mbox{div~~}a=0$. We call a pair of functions $(u,p)$ suitable weak solution to

\begin{eqnarray}\label{eq:generalizedNSE}
     \left.\begin{array}{rr}
          \partial_tu-\Delta u+a\cdot\nabla u+{\rm div~~}(a\otimes u)+u\cdot\nabla u+\nabla p&=0\\
            \mbox{div~~}u&=0
\end{array}\right\}\quad\mbox{$(x,t)\in \mathcal{O}$,}
\end{eqnarray}
if $u\in L^{\infty}_tL^2_x\cap L^2_t\dot{H}_x^1(\mathcal{O}')$ for any open subset $\mathcal{O}'\subseteq \overline{\mathcal{O}'}\Subset \mathcal{O}$, $p\in L^{3/2}_{loc}(\mathcal{O})$, such that $(u,p)$ satisfies equations (\ref{eq:generalizedNSE}) in the sense of distributions in $\mathcal{O}$, and
\begin{equation}\label{eq:generalizedEnergyInequality}
\partial_t\frac{|u|^2}{2}-\Delta \frac{|u|^2}{2}+|\nabla u|^2+\mbox{div~~}\left(\frac{|u|^2}{2}(u+a)\right)+u{\rm ~div~~}(a\otimes u)+\mbox{div~~}(up)\leq 0,
\end{equation}
in the sense of distributions. Recall that a distribution $v$ in $\mathcal{O}$ is called nonnegative if $(v,\phi)\ge 0$ for any $\phi\in C_c^{\infty}(\mathcal{O})$ with $\phi\ge 0$; $u{\rm ~div~~}(a\otimes u)$ is a distribution with
\begin{equation*}
\left(u{\rm ~div~~}(a\otimes u),\phi\right)=-\int_{\mathcal{O}}a_iu_j\partial_ju_i\phi(x,t)dxdt-\int_{\mathcal{O}}a_iu_ju_i\partial_j\phi(x,t)dxdt.
\end{equation*}
The terms in \ref{eq:generalizedEnergyInequality} make sense due to the regularity assumptions and $u\in L^{10/3}_{loc}(\mathcal{O})$ by known multiplicative inequalities.\\
The main theorem in this section can be stated as the following:
\begin{theorem}{\rm ($\epsilon$-regularity criterion)}\label{th:smallRegularity}\\
Let $(u,p)$ be a suitable weak solution to equations (\ref{eq:generalizedNSE}) in $Q_1$ with $a\in L^m(Q_1)$, $m>5$, ${\rm div~~}a=0$. Then there exists $\epsilon_0=\epsilon_0(m)>0$ with the following property: if
\begin{equation}
\left(\fint_{Q_1}~|u|^3~dxdt\right)^{1/3}+\left(\fint_{Q_1}~|p|^{3/2}~dxdt\right)^{2/3}+\left(\fint_{Q_1}~|a|^m~dxdt\right)^{1/m}\leq \epsilon_0,
\end{equation}
then $u$ is H\"{o}lder continuous in $Q_{1/2}$ with exponent $\alpha=\alpha(m)>0$ and
\begin{equation}
\|u\|_{C_{{\rm par}}^{\alpha}(Q_{1/2})}\leq C(m,\epsilon_0).
\end{equation}
\end{theorem}
\smallskip
\noindent
\textbf{Remarks:} The proof of this theorem follows the general line of presentation in \cite{EsSvSe,LaSe,Lin}. There are some additional complications due to the new terms $a\cdot\nabla u+{\rm div~~}(a\otimes u)$ as we shall see below.\\

\noindent
Before going into the proof of the theorem, we need the following two lemmas to be used below.
\begin{lemma}\label{lm:pigeonHoleInequality}
 Let $f$ be a nonnegative nondecreasing bounded function defined on $[0,1]$ with the following property:\\
for any $3/4\leq s<t<1$ and some positive constants $0<\theta<1$, $M>0$, $\beta>0$, we have
\begin{equation}
 f(s)\leq \theta f(t)+\frac{M}{(t-s)^{\beta}}.
\end{equation}
Then,
\begin{equation}
 \sup_{s\in[0,3/4]}f(s)\leq C(\theta,\beta,M),
\end{equation}
for some postive constant depending only on $\theta,~\beta,~M$.
\end{lemma}
\smallskip
\noindent
\textbf{Remarks:} The lemma is well-known, one can find a proof for example in \cite{Evans}.\\

\noindent
Our next lemma is an estimate of a generalized Stokes system.
\begin{lemma}\label{lm:estimatesForGeneralizedStokesSystem}
Let $a\in L^m(Q_1)$, with ${\rm div~~}a=0$ and $(\int_{Q_1}|a|^mdxdt)^{1/m}\leq M$, for some positive $M>0$, $m>5$, let $\lambda\in R^n$, $|\lambda|\leq M$, $f=(f_{ij})\in L^m(Q_1)$ with $(\int_{Q_1}|f|^mdxdt)^{1/m}\leq M$. Let $u\in L^{\infty}_tL^2_x\cap L^2_t\dot{H}^1_x(Q_1)$ and $p\in L^{3/2}(Q_1)$ with
\begin{equation}
 \left(\int_{Q_1}~|u|^3~dxdt\right)^{1/3}+\left(\int_{Q_1}~|p|^{3/2}~dxdt\right)^{2/3}\leq M,
\end{equation}
Assume $(u,p)$ satisfies
\begin{eqnarray}
 \left.\begin{array}{rl}
        \partial_tu-\Delta u+a\cdot\nabla u+\lambda\cdot\nabla u+{\rm div~~}(a\otimes u)+\nabla p&={\rm div~~}f\\
                                 {\rm div~~}u&=0
       \end{array}\right\}\quad{\rm in~~}Q_1,	
\end{eqnarray}
in the sense of distributions. Then $u$ is H\"{o}lder continuous in $Q_{1/2}$ with exponent $\alpha=\alpha(m)>0$ and
\begin{equation}
 \|u\|_{C_{{\rm par}}^{\alpha}(Q_{1/2})}\leq C(m,M).
\end{equation}

\end{lemma}
\smallskip
\noindent
\textbf{Proof:} The proof is based on a standard application of boostraping arguments. We sketch some of the details below. In order to use the boostraping argument, we need to show, suppose $u\in L^q(Q_R)$, $q\ge 3$, then $u\in L^{\tilde{q}}(Q_{R-\delta})$, for some $\tilde{q}>q$. Here we can assume $R>3/4$ and $\delta$ is a small positive number. Let us rewrite the equations of $(u,p)$ as
\begin{eqnarray*}
 \left.\begin{array}{rl}
        \partial_tu-\Delta u+\nabla p&={\rm div~~}(f-a\otimes u-u\otimes a-u\otimes \lambda)\\
       {\rm div~~}u&=0
       \end{array}\right\}\quad{\rm in ~~}Q_R.
\end{eqnarray*}
By H\"{o}lder inequality, we see $h:=f-a\otimes u-u\otimes a-u\otimes \lambda\in L^{\frac{mq}{m+q}}(Q_R)$. Taking divergence in the first equation, we obtain
\begin{equation*}
 \Delta p={\rm div ~~div~~}(f-a\otimes u-u\otimes a-u\otimes \lambda)\quad {\rm in~~}Q_R.
\end{equation*}
Set
\begin{equation*}
p_1=\Delta^{-1}{\rm div~~div~~}\left((f-a\otimes u-u\otimes a-u\otimes \lambda)\chi_{B_R}\right),
\end{equation*}
and write $p=p_1+p_2$. Then $\Delta p_2=0$ in $Q_R$. Recall that $R$ is in $[3/4,1]$. By elliptic estimates, we get
\begin{equation*}
 \|p_1\|_{L^{\frac{mq}{m+q}}(Q_R)}\leq C\|h\|_{L^{\frac{mq}{m+q}}(Q_R)}.
\end{equation*}
Since $\frac{mq}{m+q}>3/2$, we see $p_2$ verifies estimate
\begin{equation*}
 \|p_2\|_{L^{3/2}_tC_x^2(Q_{R-\delta/2})}\leq C(\delta,M),
\end{equation*}
with $\delta$ being a small positive number. Then,
\begin{equation*}
 \partial_tu-\Delta u=-\nabla p_1-\nabla p_2+{\rm div~~}h \quad{\rm in~~}Q_{R-\delta/2},
\end{equation*}
where $p_1,~~p_2$ and $h$ satisfy above estimates. For a smooth cutoff function $\eta$ with $\eta\equiv 1$ in $Q_{R-3\delta/4}$ and $\eta\equiv 0$ outside $Q_{R-\delta/2}$, set
\begin{eqnarray*}
v_1(\cdot,t)&=&\int_{-\infty}^te^{\Delta(t-s)}[-\nabla (p_1\eta)(\cdot,s)+{\rm div~~}(h\eta)(\cdot,s)]ds,\\
 v_2(\cdot,t)&=&\int_{-\infty}^te^{\Delta(t-s)}\nabla (p_2\eta)(\cdot,s)ds.
\end{eqnarray*}
Write $u=v_1+v_2+v_3$. By estimates of heat equation, we see $\|v_2\|_{L^{\infty}(Q_{R-\delta})}\leq C(\delta,M)$. As for $v_1$, by Young's inequality and the properties of heat kernel, we see $v_1\in L^r(Q_{R-\delta})$ for any $r>0$ such that
\begin{equation*}
 \frac{1}{r}>\frac{1}{q}+\frac{1}{m}-\frac{1}{5}.
\end{equation*}
Since $v_3$ satisfies heat equation in $Q_{R-3\delta/4}$, we see $v_3$ is smooth in $Q_{R-\delta}$. Thus in summary, we get,

\medskip
\noindent
$u\in L^q(Q_R)$ implies $u\in L^r(Q_{R-\delta})$ with $\frac{1}{r}=\frac{1}{q}-\frac{1}{2}(\frac{1}{5}-\frac{1}{m})$.

\medskip
\noindent
Since $m>5$, after applying this boostraping argument for finitely many times, we can conclude $u\in L^{r_0}(Q_{5/8})$ with $r_0$ sufficiently large such that
\begin{equation}
|a||u|\in L^{\frac{m+5}{2}}(Q_{5/8}).
\end{equation}
Then we can go back to the decompositions $v_1,~~v_2,~~$ and $v_3$, it is not difficult to verify that all of them are H\"{o}lder continuous in $Q_{1/2}$ with exponent $\alpha=\alpha(m)$. If we keep track of the constants in the above process, it's clear we also have the estimates claimed in the lemma. Alternatively, one can use the closed graph theorem with appropriate function spaces to obtain the estimates, we omit the details here. The lemma is proved.\\

\noindent
Now we can return to the proof of Theorem \ref{th:smallRegularity}. We first prove the following `oscillation lemma', which roughly speaking asserts that if $u$ is of `small oscillation' in $Q_1$, then the oscillation is even smaller in $Q_{\theta}$ for $\theta<1$.\\

\begin{lemma}{\rm (Oscillation lemma)}\label{lm:smallOscillation}\\
Let $(u,p)$ be a suitable weak solution to equations (\ref{eq:generalizedNSE}) in $Q_1$ with $a\in L^m(Q_1)$, $m>5$, ${\rm div~~}a=0$, $\|a\|_{L^m(Q_1)}\leq c$, $|(u)_{1}|\leq M$, for some small absolute number $c>0$, and some positive number $M$. Then for any $\theta\in (0,1/3)$, there exists an $\epsilon=\epsilon(\theta,M,m)>0$, $C_1(M,m)>0$, and $\alpha=\alpha(m)>0$ such that if
\begin{equation*}
Y(u,p,Q_1)+|(u)_{1}|\left(\fint_{Q_1}~|a|^m~dxdt\right)^{1/m}<\epsilon,
\end{equation*}
then
\begin{equation*}
Y(u,p,Q_{\theta})\leq C_1(M,m)\theta^{\alpha}\left(Y(u,p,Q_1)+|(u)_{1}|\left(\fint_{Q_1}~|a|^m~dxdt\right)^{1/m}\right).
\end{equation*}

\end{lemma}
\smallskip
\noindent
\textbf{Proof:} Suppose the lemma is false. Then there exists $(u_i,p_i)$ and $a_i$ with the following properties:
\begin{eqnarray*}
&&|(u_i)_{1}|\leq M, ~~\|a_i\|_{L^m(Q_1)}\leq c, ~~{\rm div~~}a_i=0,\\
&&Y(u_i,p_i,Q_1)+|(u_i)_{1}|\left(\fint_{Q_1}~|a_i|^m~dxdt\right)^{1/m}=\epsilon_i\to 0+ {\rm ~~as~~} i\to +\infty,\\
&&Y(u_i,p_i,Q_{\theta})>C_1(M,m)\theta^{\alpha}\epsilon_i,
\end{eqnarray*}
\smallskip
and $(u_i,p_i)$ satisfies equations (\ref{eq:generalizedNSE}) and inequality (\ref{eq:generalizedEnergyInequality}) with $a$ replaced by $a_i$. \\
\smallskip
Set
\begin{eqnarray*}
v_i&=&\frac{u_i-(u_i)_{1}}{\epsilon_i},\\
q_i&=&\frac{p_i-(p_i)_{1}(t)}{\epsilon_i},\\
f_i&=&\frac{a_i\otimes(u_i)_{1}}{\epsilon_i}.
\end{eqnarray*}
Then we have $(v_i)_{1}=0$, ${\rm ~div~div~}f_i=0$,
\begin{eqnarray*}
&&\left(\fint_{Q_1}|v_i|^3dxdt\right)^{1/3}+\left(\fint_{Q_1}|q_i|^{3/2}dxdt\right)^{2/3}+\left(\fint_{Q_1}|f_i|^mdxdt\right)^{1/m}\leq 1, {\rm and,}\\  &&\left(\fint_{Q_{\theta}}|v_i-(v_i)_{\theta}|^3dxdt\right)^{1/3}+\theta\left(\fint_{Q_{\theta}}|q_i-(q_i)_{\theta}(t)|^{3/2}dxdt\right)^{2/3}>C_1(M,m)\theta^{\alpha}.
\end{eqnarray*}
\smallskip
\noindent
Moreover, $(v_i,q_i)$ satisfies:
\begin{eqnarray}\label{eq:equationForScaledSolution}
         \left.\begin{array}{rr}
            \partial_tv_i-\Delta v_i+\epsilon_iv_i\cdot\nabla v_i+a_i\cdot\nabla v_i+{\rm div~~}(a_i\otimes v_i)+{\rm div~~}f_i+(u_i)_{1}\cdot\nabla v_i+\nabla q_i&=0\\
           \mbox{div~~}v_i&=0
             \end{array}\right\}
\end{eqnarray}
in the sense of distributions in $Q_1$ and
\begin{eqnarray}\label{eq:energyInequalityForScaledSolution}
\left.\begin{array}{ll}
 \partial_t\frac{|v_i|^2}{2}-\Delta\frac{|v_i|^2}{2}+|\nabla v_i|^2+\mbox{div~~}\left(\frac{|v_i|^2}{2}\left(\epsilon_i v_i+(u_i)_{1}+a_i\right)\right)&\\
+v_i{\rm ~div~~}(f_i+a_i\otimes v_i)+\mbox{div~~}v_iq_i \leq 0,&
\end{array}\right.
\end{eqnarray}
in the sense of distributions in $Q_1$. Here again the terms make sense due to our regularity assumptions and the interpretation of $v_i{\rm ~div~~}(a_i\otimes v_i+f_i)$ as the one below inequalities (\ref{eq:generalizedEnergyInequality}).\\
\smallskip
Since $v_i\in L^{\infty}_tL^2_x\cap L^2_t\dot{H}^1_x(Q_1)$ and $v_i$ satisfies equations (\ref{eq:equationForScaledSolution}), we can change the value of $v_i$ on a set of measure zero such that $t\to v_i(\cdot,t)$ is continuous from $(-1,0)$ to $L^2_w(B_1(0))$, the weak $L^2$ space.\\
From inequality (\ref{eq:energyInequalityForScaledSolution}) we obtain,
\begin{eqnarray*}
 &&\int_{B_1(0)}\frac{|v_i|^2}{2}(x,t)\phi(x,t)dx+\int_{-1}^t\int_{B_1(0)}|\nabla v_i|^2\phi(x,s)dxds\\
&&\leq \int_{-1}^t\int_{B_1(0)}\frac{|v_i|^2}{2}(\partial_t\phi+\Delta\phi)dxds+\int_{-1}^t\int_{B_1(0)}\frac{|v_i|^2}{2}[(u_i)_{1}+a_i+\epsilon_iv_i]\nabla\phi dxds\\
&&+\int_{-1}^t\int_{B_1(0)}[(f_i+a_i\otimes v_i):(\nabla v_i\phi+v_i\otimes\nabla \phi)] dxds+\int_{-1}^t\int_{B_1(0)}q_iv_i\nabla \phi dxds,
\end{eqnarray*}
for any $\phi\ge 0$ with $\phi\in C_c^{\infty}\left(B_1(0)\times(-1,t]\right)$.\\
Let us define
\begin{equation}
E_i(r)={\rm ess~}\sup_{-r^2<t\leq 0}\int_{B_r}\frac{|v_i|^2}{2}(x,t)dx+\int_{-r^2}^0\int_{B_r}|\nabla v_i|^2(x,s)dxds,
\end{equation}
for $0< r<1$. By known multiplicative inequalities we have
\begin{equation}
\|v_i\|_{L^{10/3}(Q_r)}^2\leq CE_i(r).
\end{equation}
Then for any $1/2<r_1<r_2\leq 1$, if we choose nonnegative test function $\phi$ with support in $Q_{r_2}$ appropriately, we obtain the following estimates, with the help of the above local energy estimates, H\"{o}lder estimates, and the estimates on $v_i,q_i,a_i$:
\begin{eqnarray*}
E(r_1)&\leq& \frac{C}{(r_2-r_1)^2}+\frac{C}{r_2-r_1}\int_{Q_{r_2}}\frac{|v_i|^2}{2}\left(M+|a_i|+\epsilon_i|v_i|\right)dxds\\
&+&\int_{Q_{r_2}}|f_i||\nabla v_i|+|a_i||v_i||\nabla v_i|dxdt+\frac{C}{r_2-r_1}\int_{Q_{r_2}}|f_i||v_i|+|a_i||v_i|^2+|q_i||v_i|dxdt\\
&\leq&\frac{C(M)}{(r_2-r_1)^2}+\left(\int_{Q_{r_2}}|f_i|^2dxdt\right)^{1/2}E(r_2)^{1/2}+\|a_i\|_{L^5(Q_{r_2})}\|v_i\|_{L^{10/3}(Q_{r_2})}\|\nabla v\|_{L^2(Q_{r_2})}\\
&\leq&\frac{C(M)}{(r_2-r_1)^2}+CE(r_2)^{1/2}+C\|a\|_{L^m(Q_1)}E(r_2)\\
&\leq&\frac{C(M)}{(r_2-r_1)^2}+(C\|a_i\|_{L^m(Q_1)}+1/2)E(r_2).
\end{eqnarray*}
Note that we have $\|a_i\|_{L^m(Q_1)}\leq c$ with $c$ small. So if we choose $c$ such that $Cc<1/2$, then we can apply Lemma \ref{lm:pigeonHoleInequality} and conclude that $E(3/4)\leq C(M,m)$. That is, $v_i$ are uniformly bounded in $L^{\infty}_tL^2_x\cap L^2_t\dot{H}^1_x(Q_{3/4})$. Thus by known embedding theorems and the fact that $v_i$ satisfies equations (\ref{eq:equationForScaledSolution}) (which provide crucial regularity in $t$), we can choose a subsequence of $v_i$ (which we still denote as $v_i$), such that for some $\lambda\in R^3$, $v\in L^3(Q_{3/4})$, $q\in L^{3/2}(Q_{3/4})$ and $a,f\in L^m(Q_{3/4})$ with ${\rm div~~}a=0$, we have\\
$v_i(\cdot,t)\rightharpoonup v(\cdot,t)$ weakly in $L^2(B_{3/4})$ for every $t\in (-(\frac{3}{4})^2,0)$,\\
$v_i\to v$ strongly in $L^3(Q_{3/4})$,\\
$q_i\rightharpoonup q$ weakly in $L^{3/2}(Q_{3/4})$,\\
$(u_i)_{1}\to \lambda$, $a_i\rightharpoonup a$ weakly in $L^m(Q_{3/4})$,\\
$f_i\rightharpoonup f$ in $L^m(Q_{3/4})$.\\
Moreover, we have
\begin{equation*}
\left(\int_{Q_{3/4}}|v|^3dxdt\right)^{1/3}+\left(\int_{Q_{3/4}}|q|^{3/2}dxdt\right)^{2/3}+\left(\int_{Q_{3/4}}(|f|+|a|)^mdxdt\right)^{1/m}\leq C,
\end{equation*}
and $|\lambda|\leq M$.\\
From equations (\ref{eq:equationForScaledSolution}) for $(v_i,q_i)$, we see
\begin{eqnarray}
       \left.\begin{array}{rr}
\partial_tv-\Delta v+\lambda\cdot\nabla v+{\rm div~~}(a\otimes v+v\otimes a+f)+\nabla q &=0\\
                 \mbox{div~~}v&=0
             \end{array}\right\}\quad \mbox{in $Q_{3/4}$.}
\end{eqnarray}
By Lemma \ref{lm:estimatesForGeneralizedStokesSystem} on generalized Stokes system, we see for some $\alpha=\alpha(m)>0$, $v$ is H\"{o}lder continuous in $Q_{1/2}$ with exponent $\alpha$, with respect to parabolic distance. More precisely,
\begin{equation*}
|v(x_1,t_1)-v(x_2,t_2)|\leq C(M,m)\left(|x_1-x_2|+|t_1-t_2|^{1/2}\right)^{\alpha}.
\end{equation*}
Since $v_i\to v$ strongly in $L^3(Q_{3/4})$, we see
\begin{equation*}
\left(\fint_{Q_{\theta}}|v_i-(v_i)_{\theta}|^3dxdt\right)^{1/3}\leq C(M,m)\theta^{\alpha},
\end{equation*}
for $i$ sufficiently large.\\
Note that from equations (\ref{eq:equationForScaledSolution}) we have
\begin{equation*}
-\Delta q_i=\mbox{div~div~~}\left(\epsilon_i v_i\otimes v_i+v_i\otimes a_i+a_i\otimes v_i\right).
\end{equation*}
Let $q_i=q_i^1+q_i^2$, where
\begin{equation*}
 q_i^1=(-\Delta)^{-1}\mbox{div~div~~}\left((\epsilon_i v_i\otimes v_i+v_i\otimes a_i+a_i\otimes v_i)\chi_{B_{3/4}}\right),
\end{equation*}
with $\chi_{B_{3/4}}$ being the characteristic function on $B_{3/4}$. Since $v_i$ strongly converges to $v$ in $L^3(Q_{4/3})$, we have $q_i^1-\tilde{q_i}$ strongly converges to $0$ in $L^{3/2}(Q_{4/3})$, where
\begin{equation*}
 \tilde{q_i}=(-\Delta)^{-1}\mbox{div~div~~}\left((v\otimes a_i+a_i\otimes v)\chi_{B_{3/4}}\right).
\end{equation*}
Since $v$ is bounded, we obtain by estimates of Riesz operators $\tilde{q_i}\in L^m(Q_{1/2})$. Thus
\begin{eqnarray*}
&&\theta \left(\fint_{Q_{\theta}}~|\tilde{q_i}|^{3/2}~dxdt\right)^{3/2}\\
&&\leq \theta\left(\fint_{Q_{\theta}}~|\tilde{q_i}|^m~dxdt\right)^{1/m}\leq C(M,m)\theta^{1-5/m}.
\end{eqnarray*}
Therefore, for $i$ sufficiently large, we have
\begin{equation*}
\theta \left(\fint_{Q_{\theta}}~|q_i^1|^{3/2}~dxdt\right)^{3/2}\leq C(M,m)\theta^{1-5/m}.
\end{equation*}
By definition, $\Delta q_i^2=0$ in $Q_{3/4}$ and $\left(\fint_{Q_{3/4}}~|q_i^2|^{3/2}~dxdt\right)^{2/3}\leq C$. Thus by elliptic estimates, we obtain,
\begin{eqnarray*}
 \theta\left(\fint_{Q_{\theta}}|q_i^2-(q_i^2)_{\theta}(t)|^{3/2}dxdt\right)^{2/3}&\leq& C\theta \left(\theta^{3/2}\fint_{-\theta^2}^0\|\nabla q_i^2(\cdot,t)\|_{L^{\infty}(B_{5/12})}^{3/2}dt\right)^{2/3}\\
&\leq&C\theta\left(\theta^{-1/2}\int_{-\theta^2}^0\int_{B_{1/2}(0)}|q|^{3/2}dxdt\right)^{2/3}\\
&\leq&C\theta^{2/3}.
\end{eqnarray*}
Therefore, summarizing the above, we see
\begin{equation*}
\theta\left(\fint_{Q_{\theta}}|q_i-(q_i)_{\theta}(t)|^{3/2}dxdt\right)^{2/3}\leq C(M,m)\theta^{\min\{2/3,1-5/m\}},
\end{equation*}
for $i$ sufficiently large. This, together with the estimates on $v_i$, shows
\begin{equation*}
Y(v_i,q_i,Q_{\theta})\leq C(M,m)\theta^{\alpha},
\end{equation*}
 for $i$ sufficiently large, if we choose $\alpha(m)$ sufficiently small. This contradicts $Y(v_i,q_i,Q_{\theta})\ge C_1(M,m)\theta^{\alpha}$ if we choose $C_1(M,m)>2C(M,m)$. Thus the lemma is proved.\\

The above lemma admits the following iterations.
\begin{lemma}{\rm (Iteration of the oscillation lemma)}\label{lm:iterationOfOscillation}\\
Let $(u,p)$, $M$, $\epsilon(\theta,M,m)$, $C_1(M,m)$, $\alpha(m)$, $c$ and $a$ be as in the above lemma, with $|(u)_{Q_1}|\leq M/2$. Let $\beta=\alpha/2$. Choose $\theta\in (0,1/3)$ such that $C_1(M,m)\theta^{\alpha-\beta}<1$, and $\theta<c_1$ with $c_1=c_1(M,m)$ being some small number. Then there exists $\epsilon_{\ast}(\theta,M,m)$ sufficiently small, such that if
\begin{equation}
Y(u,p,Q_1)+M\left(\fint_{Q_1}|a|^mdxdt\right)^{1/m}<\epsilon_{\ast},                                                                                 \end{equation}
then for any $k=1,2,\dots$, we have
\begin{eqnarray}
&&|(u)_{Q_{\theta^{k-1}}}|\leq M,\\
&&Y(u,p,Q_{\theta^{k-1}})+|(u)_{\theta^{k-1}}|\left(\fint_{Q_{\theta^{k-1}}}|a|^mdxdt\right)^{1/m}\theta^{k-1}<\epsilon_{\ast}\leq\epsilon(\theta,M,m), \\
&&Y(u,p,Q_{\theta^k})\leq\theta^{\beta}\left(Y(u,p,Q_{\theta^{k-1}})+|(u)_{\theta^{k-1}}|\left(\fint_{Q_{\theta^{k-1}}}|a|^mdxdt\right)^{1/m}\theta^{k-1}\right).
\end{eqnarray}
\end{lemma}
\noindent
\textbf{Proof:} We prove the lemma by induction.\\
For $k=1$, the conclusion follows from Lemma \ref{lm:smallOscillation}, if we choose $\epsilon_{\ast}$ such that $\epsilon_{\ast}<\epsilon(\theta,M,p)$. Suppose the conclusion is true for $k\leq k_0$, $k_0\ge 1$, we show it remains true for $k=k_0+1$.\\
By induction
\begin{eqnarray*}
&&|(u)_{Q_{\theta^{k-1}}}|\leq M,\\
&&Y(u,p,Q_{\theta^{k-1}})+|(u)_{\theta^{k-1}}|\left(\fint_{Q_{\theta^{k-1}}}|a|^mdxdt\right)^{1/m}\theta^{k-1}<\epsilon_{\ast},\\
&&Y(u,p,Q_{\theta^{k}})\leq \theta^{\beta}\left(Y(u,p,Q_{\theta^{k-1}})+|(u)_{\theta^{k-1}}|\left(\fint_{Q_{\theta^{k-1}}}|a|^mdxdt\right)^{1/m}\theta^{k-1}\right)\leq\theta^{\beta}\epsilon_{\ast},
\end{eqnarray*}
for all $k\leq k_0$. Thus,
\begin{eqnarray*}
 Y(u,p,Q_{\theta^k})&\leq &\theta^{\beta}\left(Y(u,p,Q_{\theta^{k-1}})+\theta^{k-1}M\left(\fint_{Q_{\theta^{k-1}}}|a|^mdxdt\right)^{1/m}\right)\\
&\leq &\theta^{\beta}\left(Y(u,p,Q_{\theta^{k-1}})+\theta^{(k-1)(1-5/m)}M\left(\fint_{Q_1}|a|^mdxdt\right)^{1/m}\right)\\
&\leq&\theta^{\beta}Y(u,p,Q_{\theta^{k-1}})+\theta^{k\beta_1}\epsilon_{\ast}
\end{eqnarray*}
for all $k\leq k_0$, with $\beta_1=\min\{\beta,1-5/m\}$. Simple calculations with a repeated use of the above inequalities show:
\begin{equation*}
 Y(u,p,Q_{\theta^k})\leq \theta^{k\beta}Y(u,p,Q_1)+k\theta^{k\beta_1}\epsilon_{\ast}, \quad\forall k\leq k_0.
\end{equation*}
Thus,
\begin{eqnarray*}
|(u)_{Q_{\theta^{k_0}}}|&\leq &\sum_{k=1}^{k_0}|(u)_{Q_{\theta^k}}-(u)_{Q_{\theta^{k-1}}}|+|(u)_{Q_1}|\\
&\leq&\sum_{k=1}^{k_0}\left(\fint_{Q_{\theta^k}}|u-(u)_{Q_{\theta^{k-1}}}|^3dxdt\right)^{1/3}+|(u)_{Q_1}|\\
&\leq&\theta^{-5/3}\sum_{k=1}^{k_0}Y(u,p,Q_{\theta^{k-1}})+|(u)_{Q_1}|\\
&\leq&\theta^{-5/3}\sum_{k=1}^{k_0}\left(\theta^{(k-1)\beta}\epsilon_{\ast}+\epsilon_{\ast}(k-1)\theta^{(k-1)\beta_1}\right)+M/2\\
&\leq&\theta^{-5/3}(1-\theta^{\beta})^{-1}\epsilon_{\ast}+\theta^{-5/3}\epsilon_{\ast}C(\beta_1,\theta)+M/2.
\end{eqnarray*}
If we choose $\epsilon_{\ast}=\epsilon_{\ast}(\theta,M,m)$ to be sufficiently small, we see
\begin{equation*}
 |(u)_{Q_{\theta^{k_0}}}|\leq M.
\end{equation*}

Moreover,
\begin{eqnarray*}
&&Y(u,p,Q_{\theta^{k_0}})+\theta^{k_0}|(u)_{\theta^{k_0}}|\left(\fint_{Q_{\theta^{k_0}}}|a|^mdxdt\right)^{1/m}\\
&&\leq \theta^{\beta}\epsilon_{\ast}+\theta^{(1-5/m)k_0}\epsilon_{\ast}<\epsilon_{\ast},
\end{eqnarray*}
if we choose $\theta<c(M,m)$ to be sufficiently small. Set
\begin{eqnarray*}
&&u(x,t)=\frac{1}{\theta^{k_0}}v(\frac{x}{\theta^{k_0}},\frac{t}{\theta^{2k_0}}),\\ &&p(x,t)=\frac{1}{\theta^{2k_0}}q(\frac{x}{\theta^{k_0}},\frac{t}{\theta^{2k_0}}), {\rm ~~and}\\ &&a(x,t)=\frac{1}{\theta^{k_0}}b(\frac{x}{\theta^{k_0}},\frac{t}{\theta^{2k_0}}).
 \end{eqnarray*}
One can verify that $(v,q)$ is a suitable weak solution to equations (\ref{eq:generalizedNSE}) with $a$ replaced by $b$ in $Q_1$. Moreover,
\begin{eqnarray*}
 &&Y(v,q,Q_1)+|(v)_{Q_1}|\left(\fint_{Q_1}|b|^mdxdt\right)^{1/m}\\
&&=\theta^{k_0}\left(Y(u,p,Q_{\theta^{k_0}})+\theta^{k_0}|(u)_{Q_{\theta^{k_0}}}|\left(\fint_{Q_{\theta^{k_0}}}|a|^mdxdt\right)^{1/m}\right)<\epsilon_{\ast},\\
&&\left(\int_{Q_1}|b|^mdxdt\right)^{1/m}\leq\theta^{k_0-\frac{5k_0}{m}}\left(\int_{Q_1}|a|^mdxdt\right)^{1/m}<c.
\end{eqnarray*}
Thus, by Lemma \ref{lm:smallOscillation}, we obtain,
\begin{equation}
Y(v,q,Q_{\theta})\leq \theta^{\beta}\left(Y(v,q,Q_1)+|(v)_{Q_1}|\left(\fint_{Q_1}|b|^mdxdt\right)^{1/m}\right),
\end{equation}
that is, 
\begin{equation}
Y(u,p,Q_{\theta^{k_0+1}})\leq\theta^{\beta}\left(Y(u,p,Q_{\theta^{k_0}})+|(u)_{\theta^{k_0}}|\left(\fint_{Q_{\theta^{k_0}}}|a|^mdxdt\right)^{1/m}\theta^{k_0}\right). 
\end{equation}
The lemma is then proved.\\

By translation and dilation, we obtain the following corollary.
\begin{corollary}\label{cor:rescaledOscillation}
Let $(u,p)$ be a suitable weak solution to equations (\ref{eq:generalizedNSE}) in $Q(R,z_0)$, with $a\in L^m(Q(R,z_0))$, ${\rm div~~}a=0$, $|(u)_{Q(R,z_0)}|R<M/2$, $\theta$, $\beta$ are as in the above. Then there exists $\epsilon_{\ast}=\epsilon_{\ast}(\theta,M,m)$ such that
\begin{equation*}
RY(u,p,Q(R,z_0))+RM\left(\fint_{Q(R,z_0)}|a|^mdxdt\right)^{1/m}<\epsilon_{\ast}
\end{equation*}
 implies, for $k\ge 1$ {\rm :}
\smallskip
\begin{eqnarray*}
&&R|(u)_{Q(\theta^{k-1}R,z_0)}|\leq M, {\rm ~~and}\\
&&Y(u,p,Q(\theta^{k}R,z_0))\\
&&\leq\theta^{\beta}\left(Y\left(u,p,Q(\theta^{k-1}R,z_0)\right)+R\theta^{k-1}|(u)_{Q_{\theta^{k-1}R}}|\left(\fint_{Q(\theta^{k-1}R,z_0)}|a|^mdxdt\right)^{1/m}\right).
\end{eqnarray*}
 \end{corollary}
\smallskip
\noindent
\textbf{Proof of Theorem \ref{th:smallRegularity}:}\\
It is clear if we choose $\epsilon_0$ sufficiently small, we can apply Corollary \ref{cor:rescaledOscillation} in $Q(1/2,z_0)$ for any $z_0\in Q_{1/2}$. Note that $|(u)_{Q_{\theta^{k}R}}|$ is bounded and $m>5$. Thus we can conclude
\begin{equation*}
 Y(u,p,Q(z_0,Q_{\theta^k}))\leq C(\theta,M,m)\theta^{k\alpha},
\end{equation*}
for some $\alpha=\alpha(m)$, where we can choose $M<1$, $\theta=\theta(M,m)=\theta(m)$. (There is a slight abuse of notation, in particular, this $\alpha$ is smaller than those appearing in the oscillation lemma). In particular,
\begin{equation*}
  \left(\fint_{Q(\theta^k,z_0)}|u-(u)_{Q(\theta^k,z_0)}|^3dxdt\right)^{1/3}\leq C(\theta,M,m)\theta^{k\alpha},
\end{equation*}
for all $z_0\in Q_{1/2}$ and $k\ge 1$. By Campanato's lemma, we conclude $u$ is H\"{o}lder continuous in $Q_{1/2}$. The theorem is proved.\\

In applications, it is cumbersome to have the ``smallness condition'' on $a$. We can remove this condition and get the following theorem.
\begin{theorem}{\rm (Improved $\epsilon$-regularity criteria)}\label{th:improvedSmallRegularity}\\
Let $(u,p)$ be a suitable weak solution to equations (\ref{eq:generalizedNSE}) in $Q_1$, with $a\in L^m(Q_1)$, ${\rm div~~}a=0$, $\|a\|_{L^m(Q_1)}\leq M$, for some $M>0$ and $m>5$. Then there exists $\epsilon_1=\epsilon_1(m,M)>0$ with the following properties: if
\begin{equation*}
\left(\fint_{Q_1}|u|^3dxdt\right)^{1/3}+\left(\fint_{Q_1}|p|^{3/2}dxdt\right)^{2/3}\leq \epsilon_1,
\end{equation*}
then $u$ is H\"{o}lder continuous in $Q_{1/2}$ with exponent $\alpha=\alpha(m)>0$ and
\begin{equation}
\|u\|_{C_{{\rm par}}^{\alpha}(Q_{1/2})}\leq C(m,\epsilon_1,M)=C(m,M).
\end{equation}
\end{theorem}
\smallskip
\noindent
\textbf{Proof:} Choose $0<R_0<1/2$, a small positive number to be determined below. For any $z_0=(x_0,t_0)\in Q_{1/2}$, we would like to apply a scaled version of Theorem \ref{th:smallRegularity} for $(u,p)$ in $Q(R,z_0)$. Set
\begin{eqnarray*}
&&u(x,t)=\frac{1}{R_0}v(\frac{x-x_0}{R_0},\frac{t-t_0}{R_0^2}), \\
&&p(x,t)=\frac{1}{R_0^2}q(\frac{x-x_0}{R_0},\frac{t-t_0}{R_0^2}),\\
&&a(x,t)=\frac{1}{R_0}b(\frac{x-x_0}{R_0},\frac{t-t_0}{R_0^2}).
\end{eqnarray*}
We see that $(v,q)$ is a suitable weak solution to equations (\ref{eq:generalizedNSE}) with $a$ replaced by $b$ in $Q_1$. Moreover,
\begin{equation*}
\|b\|_{L^m(Q_1)}\leq R_0^{1-5/m}\|a\|_{L^m(Q(R_0,z_0))}\leq CR_0^{1-5/m}M,
\end{equation*}
and
\begin{eqnarray*}
&&\left(\fint_{Q_1}|v|^3dxdt\right)^{1/3}+\left(\fint_{Q_1}|q|^{3/2}dxdt\right)^{2/3}\\
&&=R_0\left(\fint_{Q(R_0,z_0)}|u|^3dxdt\right)^{1/3}+\left(\fint_{Q(R_0,z_0)}|p|^{3/2}dxdt\right)^{2/3}R_0^2\\
&&\leq C (R_0R_0^{-5/3}+R_0^2R_0^{-10/3})\epsilon_1\leq CR_0^{-4/3}\epsilon_1.
\end{eqnarray*}
Thus,
\begin{eqnarray*}
&&\left(\fint_{Q_1}|v|^3dxdt\right)^{1/3}+\left(\fint_{Q_1}|q|^{3/2}dxdt\right)^{2/3}+\left(\fint_{Q_1}|b|^mdxdt\right)^{1/m}\\
&&\leq R_0^{1-5/m}M+CR_0^{-4/3}\epsilon_1.
\end{eqnarray*}
Thus, if we choose $R_0$ such that $R_0^{1-5/m}M<\epsilon_0/2$, fix $R_0$, $R_0=R_0(M,m)$ and choose $\epsilon_1$ such that $CR_0^{-4/3}\epsilon_1<\frac{\epsilon_0}{2}$. Then we can apply Theorem \ref{th:smallRegularity} to $(v,q)$ and conclude $v$ is H\"{o}lder continuous in $Q_{1/2}$. Scale back and collect all constants, the theorem is then proved.\\

\end{section}

\begin{section}{Local in space near initial time smoothness of Leray solutions}
In this section, we use the `$\epsilon$-regularity' theorem proved in the last section to study the local in space near initial time smoothness of the so called Leray solutions. Our setting is as follows.\\
Let $u_0\in L^2_{loc}(R^3)$ with $\mbox{div~~}u_0=0$ and $\sup_{x_0\in R^3}\int_{B_1(x_0)}|u_0|^2dx<\infty$. We recall the definition of Leray solutions in \cite{LemRie}, see also \cite{HaSv}.
\begin{definition}{(Leray solution)} A vector field
$u\in L^2_{loc}(R^3\times [0,\infty))$ is called a Leray solution to Navier-Stokes equations with initial data $u_0$ if it satisfies:

\smallskip
\noindent
i)   ${\rm ess}\sup_{0\leq t<R^2}\sup_{x_0\in R^3}\int_{B_R(x_0)}\frac{|u|^2}{2}(x,t)dx+\sup_{x_0\in R^3}\int_0^{R^2}\int_{B_R(x_0)}|\nabla u|^2dxdt<\infty$, and
$$\lim_{|x_0|\to \infty}\int_0^{R^2}\int_{B_R(x_0)}|u|^2(x,t)dxdt=0,$$
for any $R<\infty$.

\smallskip
\noindent
ii)  for some distribution $p$ in $R^3\times (0,\infty)$, $(u,p)$ verifies Navier Stokes equations
\begin{eqnarray}\label{eq:NavierStokesEquation}
 \left.\begin{array}{rl}
        \partial_tu-\Delta u+u\cdot\nabla u+\nabla p&=0\\
         \mbox{div~~}u&=0
       \end{array}\right\}\quad\mbox{in $R^3\times(0,\infty)$,}
\end{eqnarray}
in the sense of distributions and for any compact set $K\subseteq R^3$, $\lim_{t\to 0+}\|u(\cdot,t)-u_0\|_{L^2(K)}=0$.

\smallskip
\noindent
iii)  $u$ is suitable in the sense of Caffarelli-Kohn-Nirenberg, more precisely, the following local energy inequality holds:
\begin{equation}\label{eq:no2}
\int_0^{\infty}\int_{R^3}|\nabla u|^2\phi(x,t)dxdt\leq \int_0^{\infty}\int_{R^3}\frac{|u|^2}{2}(\partial_t\phi+\Delta \phi)+\frac{|u|^2}{2}u\cdot \nabla \phi+pu\cdot\nabla\phi dxdt
 \end{equation}
for any smooth $\phi\ge 0$ with ${\rm ~supp~~} \phi\Subset R^3\times (0,\infty)$. The set of all Leray solutions starting from $u_0$ will be denoted as $\mathcal{N}(u_0)$.
\end{definition}
\smallskip
\noindent
\textbf{Remarks:} For general existence result of Leray solutions, see \cite{Cald,KikuSere,LemRie}. For us, the a priori estimates of Leray solutions below are more important, since in our situation when $u_0$ is usually better than that in \cite{LemRie}, the existence can be proved in simpler ways. In the case the initial data is in $L^2(R^3)$, the notion of Leray-Hopf weak solutions is often used (see \cite{LaSe} for example). The difference is that Leray-Hopf weak solutions belong to $L^{\infty}_tL^2_x\cap L^2_t\dot{H}^{1}_x(R^3\times [0,\infty))$. It is clear that our definition includes such solutions. Note that we impose a decay condition on $u$ in i). This condition allows us to calculate $p$ in the following way: $\forall B_r(x_0)\times (0,t_{\ast})\subseteq R^3\times(0,\infty)$, take a smooth cutoff function $\phi$ with $\phi|_{B_{2r}(x_0)}=1$,  then there exists a function $p(t)$ depending only on $x_0,r,t,\phi$ (we suppress the dependence on $x_0,r,\phi$ in our notation) such that for $(x,t)\in B_r(x_0)\times(0,t_{\ast})$
\begin{equation}\label{eq:no3}
 p(x,t)=-\Delta^{-1}\mbox{div~div}(u\otimes u\phi)-\int_{R^3}\left(k(x-y)-k(x_0-y)\right)u\otimes u(y,t)\left(1-\phi(y)\right)dy+p(t)
\end{equation}
where $k(x)$ is the kernel of $\Delta^{-1}\mbox{div~div}$. \\
The right hand side is well defined since $u$ satisfies the estimates in i) and
\begin{equation}
|k(x-y)-k(x_0-y)|=O(\frac{1}{|x_0-y|^4}) {\rm ~~as~~} |y|\to \infty.
\end{equation}
The situation is similar to extending the domain of singular integrals to bounded functions, see for example \cite{LemRie} and \cite{Stein}.\\

For Leray solution $u\in \mathcal{N}(u_0)$, we have the following a priori estimates, first proved in \cite{LemRie}, see also a simpler proof in \cite{HaSv}. These estimates have played an important role in \cite{HaSv,RuSv}, see also \cite{Seregin}.\\
\begin{lemma}{\rm (A priori estimate for Leray solutions)}\label{lm:aprioriEstimatesForLeraySolution}\\
Let $\alpha=\sup_{x_0\in R^3}\int_{B_R(x_0)}\frac{|u_0|^2}{2}(x)dx<\infty$ for some $R>0$ and let $u$ be a Leray solution with initial data $u_0$. Then there exists some small absolute number $c>0$ such that for $\lambda$ satisfying $0<\lambda\leq c\min\{\alpha^{-2}R^2,1\}$, we have
\begin{equation}\label{eq:no4}
 {\rm ess}\sup_{0\leq t\leq \lambda R^2}\sup_{x_0\in R^3}\int_{B_R(x_0)}\frac{|u|^2}{2}(x,t)dx+\sup_{x_0\in R^3}\int_0^{\lambda R^2}\int_{B_R(x_0)}|\nabla u|^2(x,t)dxdt\leq C\alpha.
\end{equation}
\end{lemma}
\smallskip
\noindent
\textbf{Remarks:} Note that from the formula (\ref{eq:no3}) and the a priori estimate of $u$, we get the following estimate for $p$ which will be useful:
\begin{equation}
\sup_{x_0\in R^3}\int_0^{\lambda R^2}\int_{B_R(x_0)}|p-p(t)|^{3/2}dxdt\leq C \alpha^{3/2}R^{1/2}.
\end{equation}
In the above estimate on $p$, more precisely, $p(t)=p_{x_0,R}(t)$. That is, we need to choose some appropriate constants $p_{x_0,R}(t)$ to satisfy the inequality. The point here is that such constants depending on $x_0,R,t$ exist. This remark is effective throughout the paper.\\

Now we can prove our our first important result.\\
\begin{theorem}\label{th:localRegularityForLeraySolution}
Let $u_0\in L^2_{loc}(R^3)$ with $\sup_{x_0\in R^3}\int_{B_1(x_0)}|u|^2(x)dx\leq \alpha<\infty$. Suppose $u_0$ is in $L^m(B_2(0))$ with $\|u_0\|_{L^m(B_2(0))}\leq M<\infty$ and $m>3$.  Let us decompose\footnote{Such decomposition is well-known. One can for example first localize $u_0$ using a smooth cutoff function, and then use Bogovskii's lemma to deal with the divergence-free condition. See for example \cite{Bogovski, GeiHechH}.} $u_0=u_0^1+u_0^2$ with ${\rm div~~}u_0^1=0$, $u_0^1|_{B_{4/3}}=u_0$, ${\rm ~supp~~} u_0^1\Subset B_2(0)$ and $\|u_0^1\|_{L^m(R^3)}\leq C(M,m)$.  Let $a$ be the locally in time defined mild solution to Navier-Stokes equations with initial data $u_0^1$.  Then there exists a positive $T=T(\alpha,m,M)>0$, such that any Leray solution $u\in \mathcal{N}(u_0)$ satisfies:\\
$u-a\in C_{{\rm par}}^{\gamma}(\overline{B_{1/2}}\times[0,T])$, and $\|u-a\|_{C_{{\rm par}}^{\gamma}(\overline{B_{1/2}}\times[0,T])}\leq C(M,m,\alpha)$, for some $\gamma=\gamma(m)\in (0,1)$.
\end{theorem}
\smallskip
\noindent
\textbf{Remark:} We can certainly choose $T(M)>0$ such that $a$ is defined on $R^3\times[0,T(M)]$. The point of the theorem is that regularity of solution to Navier-Stokes equations depends locally on initial data, as least when H\"{o}lder continuity is concerned. \\
\noindent
\textbf{Proof:} By assumption $a$ solves the Cauchy problem for Navier-Stokes equations with initial data $u_0^1$ in $R^3\times[0,T_1]$, where $T_1=T_1(M,m)$, namely:
\begin{eqnarray}
 \left.\begin{array}{rl}
        \partial_ta-\Delta a+a\cdot\nabla a+\nabla\tilde{p}&=0\\
         \mbox{div~~}a&=0
       \end{array}\right\}&&\quad\mbox{in $R^3\times(0,T_1)$,}\\
\mbox{and}\quad a(\cdot,0)=u_0^1&&.
\end{eqnarray}
It is well-known how to construct the so called mild solution to Navier-Stokes equations, see for example \cite{Kato,KoTa,Leray}. In our case, it is even simpler, since $u_0^1\in L^m$ with $m>3$ is subcrical with respect to the natural scaling of the equation. We can follow the arguments in the Appendix of \cite{EsSvSe}, and obtain $a\in L^{\frac{5m}{3}}(R^3\times(0,T_1)))$ with $\|a\|_{L^{\frac{5m}{3}}(R^3\times(0,T_1))}\leq C M$. Note that $\frac{5m}{3}>5$ since $m>3$. Moreover, by the estimates on $a$ and by treating the nonlinear term as pertubation, we can recover a local energy estimate for $a$:
\begin{equation*}
 {\rm ess~~}\sup_{0<t<T_1}\int_{B_1(x_0)}\frac{|a|^2}{2}(x,t)dx+\int_{B_1(x_0)}\int_0^{T_1}|\nabla a|^2(x,t)dxdt\leq C(M,m),
\end{equation*}
for any $x_0\in R^3$.
Write $u=a+v$, we can verify that $v$ satisfies:
\begin{eqnarray}
 \left.\begin{array}{rl}
        \partial_tv-\Delta v+v\cdot\nabla v+a\cdot\nabla v+{\rm ~div~~}(a\otimes v)+\nabla q&=0\\
       \mbox{div~~}v&=0
       \end{array}\right\}
\end{eqnarray}
in the sense of distributions in $R^3\times(0,T_1)$, here $q=p-\tilde{p}$ with $p$ being the associated pressure for $u$; and the local energy inequality
\begin{equation*}
\partial_t\frac{|v|^2}{2}-\Delta\frac{|v|^2}{2}+|\nabla v|^2+\mbox{div~}(\frac{|v|^2}{2}(v+a))+v{\rm ~div~~}(a\otimes v)+\mbox{div~~}(vq)\leq 0,
\end{equation*}
in the sense of distributions in $R^3\times(0,T_1)$;
\begin{equation*}
\lim_{t\to 0+}\|v(\cdot,t)-u_0^2\|_{L^2(B_1(x_0))}=0, {\rm ~~for ~~any~~} x_0\in R^3.
\end{equation*}
Note also that $u_0^2|_{B_{4/3}}\equiv 0$. Since $(u,p)$ satisfies the a priori estimates in Lemma \ref{lm:aprioriEstimatesForLeraySolution} (and the remarks below it), $(a,\tilde{p})$ is regular, we obtain the following estimates for $(v,q)$ in $B_2(0)\times[0,T_2)$, $T_2=T_2(\alpha,M,m)$:
\begin{eqnarray*}
 &&{\rm ess}\sup_{0<t<T_2}\frac{1}{2}\int_{B_2(0)}|v|^2(x,t)dx+\int_0^{T_2}\int_{B_2(0)}|\nabla v|^2(x,s)dxds\\
&&+\left(\int_0^{T_2}\int_{B_2(0)}|q|^{3/2}dxds\right)^{2/3}\leq C(\alpha,m,M).
\end{eqnarray*}
From the local energy inequality for $v$, and $\lim_{t\to 0+}\|v(\cdot,t)\|_{L^2(B_{4/3}(0))}=0$, we obtain
\begin{eqnarray*}
&&\frac{1}{2}\int_{B_{4/3}}|v|^2(x,t)\phi(x)dx+\int_0^t\int_{B_{4/3}}|\nabla v|^2(x,t)\phi(x)dxds\\
&&\leq\int_0^t\int_{B_{4/3}}\frac{|v|^2}{2}\Delta \phi dxds+\int_0^t\int_{B_{4/3}}\frac{|v|^2}{2}(v+a)\nabla \phi dxds\\
&&+\int_0^t\int_{B_{4/3}}[a\otimes v:(\nabla v \phi+v\otimes\nabla\phi)]+qv\cdot\nabla\phi dxds,
\end{eqnarray*}
where $\phi\in C_c^{\infty}(B_{4/3})$, $\phi|_{B_1}\equiv 1$, $\phi\ge 0$.\\
By multiplicative inequalities, we know
\begin{equation*}
\left(\int_0^{T_2}\int_{B_2(0)}|v|^{10/3}dxdt\right)^{3/10}\leq C(\alpha,m,M).
\end{equation*}
Thus from the above, we see by Schwartz inequality:
\begin{equation*}
\frac{1}{2}\int_{B_1(0)}|v|^2(x,t)dx+\int_0^t\int_{B_1(0)}|\nabla v|^2(x,s)dxds\leq C(\alpha,m,M)t^{\min\{1/30,\frac{m-3}{5m}\}},
\end{equation*}
for $t<T_2$. From
\begin{equation*}
\Delta q=-\mbox{div~div~}(v\otimes v+a\otimes v+v\otimes a),
\end{equation*}
we can see $q\in L^{5/3}_{loc}$. Thus
\begin{equation*}
\left(\int_0^t\int_{B_1(0)}|q|^{3/2}dxds\right)^{2/3}\leq C(\alpha,m,M)t^{1/15}.
\end{equation*}
The importance of these estimates lies in the fact that they provide crucial ``quantitative'' information on the decay in time as $t\to 0+$. Now for $t_0$ fixed, whose precise value is to be determined later, extend $v,q$ to $B_1(0)\times(-1+t_0,t_0]$ by setting $v=0$, $q=0$ for $(x,t)\in B_1\times(-1+t_0,0]$. Extend $a$ to $B_1(0)\times(-1+t_0,t_0]$ by setting $a(t,x)=0$ for $t<0$. The extended function $(v,q)$ is a suitable weak solution to the generalized Navier-Stokes equations (\ref{eq:generalizedNSE}) with the extended $a$ in $B_1(0)\times[-1+t_0,t_0]$. Note here that
\begin{equation*}
 \lim_{t\to 0+}\|v(\cdot,t)\|_{L^2(B_1(0))}=0
\end{equation*}
plays a crucial role: it guarantees that $\partial_tv$ and $\partial_t\frac{|v|^2}{2}$ will not cause any problem across $t=0$.
Then clearly if we choose $t_0=t_0(\alpha,m,M)$ sufficiently small, we can apply Theorem \ref{th:improvedSmallRegularity} and conclude $v$ is H\"{o}lder continuous in $B_{1/2}\times[0,t_0]$, with $\|v\|_{C_{{\rm par}}^{\gamma}(B_{1/2}\times[0,t_0])}\leq C(\alpha,m,M)$, for some $\gamma=\gamma(m)$. The theorem is proved.\\

For applications below, we state the following simple (and certainly well-known) lemma for heat equation without proof.
\begin{lemma}\label{lm:holderForHeatEquation}
 We have the following estimates:\\
1. If $u_0\in C^{\beta}(R^3)$ for some $\beta\in(0,1)$, then $e^{\Delta t}u_0(x)\in C_{{\rm par}}^{\beta}(R^3\times[0,1])$, with
\begin{equation}
 \|e^{\Delta t}u_0(x)\|_{C_{{\rm par}}^{\beta}(R^3\times[0,1])}\leq C\|u_0\|_{C^{\beta}(R^3)}.
\end{equation}
2. If $f\in L^{\infty}(R^3\times[0,1])$, then $\int_0^t\nabla e^{\Delta(t-s)}f(\cdot,s)ds\in C_{{\rm par}}^{\beta}(R^3\times[0,1])$ for any $\beta\in(0,1)$, and
\begin{equation}
 \|\int_0^t\nabla e^{\Delta(t-s)} f(\cdot,s)ds\|_{C_{{\rm par}}^{\beta}(R^3\times[0,1])}\leq C(\beta)\|f\|_{L^{\infty}(R^3\times[0,1])}.
\end{equation}

\end{lemma}

The above theorem implies the following result.
\begin{theorem}{\rm (Local H\"{o}lder regularity of Leray solutions)}\label{th:localHolderLeraySolution}\\
Let $u_0\in L^2_{loc}(R^3)$ with $\sup_{x_0\in R^3}\int_{B_1(x_0)}|u|^2(x)dx\leq \alpha<\infty$. Suppose $u_0$ is in $C^{\gamma}(B_2(0))$ with $\|u_0\|_{C^{\gamma}(B_2(0))}\leq M<\infty$. Then there exists a positive $T=T(\alpha,\gamma,M)>0$, such that any Leray solution $u\in \mathcal{N}(u_0)$ satisfies:
\begin{equation}
u\in C_{{\rm par}}^{\gamma}(\overline{B_{1/4}}\times[0,T]), {\rm ~~and~~} \|u\|_{C_{{\rm par}}^{\gamma}(\overline{B_{1/4}}\times[0,T])}\leq C(M,\alpha,\gamma).
\end{equation}
\end{theorem}
\smallskip
\noindent
\textbf{Proof:} Let us decompose $u_0=u_0^1+u_0^2$ with ${\rm div~~}u_0^1=0$, $u_0^1|_{B_{4/3}(0)}=u_0$, ${\rm ~supp~~} u_0^1\Subset B_2(0)$ and $\|u_0^1\|_{C^{\gamma}(R^3)}\leq CM$. Let $a$ be the mild solution to Navier-Stokes equations with initial data $u_0^1$ in $R^3\times(0,T(M))$. Then Theorem \ref{th:localRegularityForLeraySolution} implies that $u-a$ is H\"{o}lder continuous with some exponent $\beta\in(0,\gamma)$ in $B_{1/2}\times[0,T]$ with some $T=T(\alpha,\gamma,M)\in (0,T(M))$. Since the initial data $u_0^1$ for $a$ is in $C^{\gamma}$, it is not difficult to show that $a\in C_{{\rm par}}^{\gamma}(R^3\times(0,T))$. Thus $u$ is H\"{o}lder continuous with exponent $\beta$ in $B_{1/2}\times[0,T(M)]$ . By using a routine boostraping argument, one can improve the exponent to $\gamma$. Since this argument will be used one more time below, we sketch some of the details here for the reader's convenience. Note that $u$ is H\"{o}lder continuous in $\overline{B_{1/2}}\times[0,T]$, thus from the representation formula (\ref{eq:no3}) for $p$ and estimates for Riesz transform, we know $p$ is bounded in $\overline{B_{7/16}}\times[0,T]$ modulo some function $p(t)$. Now rewrite the equation for $u$ as
\begin{equation}
 \partial_tu-\Delta u=-{\rm div~~}(u\otimes u)-\nabla p.
\end{equation}
Choose a smooth cutoff function $\eta$ with $\eta\equiv 1$ on $\overline{B_{3/8}}$ and $\eta\equiv 0$ outside $B_{7/16}$. Write
\begin{eqnarray*}
 u_1(\cdot,t)&=&\int_0^te^{\Delta (t-s)}[-{\rm div~~}(u\otimes u\eta)-\nabla (p\eta)](\cdot,s)ds,\\
u_2(\cdot,t)&=&e^{\Delta t}(u_0\eta).
\end{eqnarray*}
Let $u=u_1+u_2+u_3$. By Lemma \ref{lm:holderForHeatEquation} we see $u_1$ and $u_2$ are H\"{o}lder continuous with exponent $\beta$. Note that $u_3$ satisfies
\begin{equation*}
 \partial_tu_3-\Delta u_3=0\quad {\rm in~~}B_{3/8}\times[0,T],
\end{equation*}
and $u_3(\cdot,0)|_{B_{3/8}}=0$. Thus $u_3$ is smooth in $\overline{B_{1/4}}\times[0,T]$. In summary $u$ is H\"{o}lder continuous in $\overline{B_{1/4}}\times[0,T]$ with exponent $\beta$. Then the theorem is proved.
\end{section}

\begin{section}{Estimates of forward self similar solutions to Navier-Stokes and Stokes equations}
In this section, we start to study forward self similar solutions to Navier-Stokes equations and a related nonhomogeneous Stokes system. Our setting is as follows.\\

Let $u$ be a Leray solution with initial data $u_0$. Suppose $\lambda u_0(\lambda x)=u_0(x)$, $\lambda u(\lambda x,\lambda^2 t)=u(x,t)$ for any $\lambda>0$. We also assume $u_0|_{\partial B_1(0)}\in C^{\infty}(\partial B_1(0))$. Then it is easy to see
\begin{equation*}
|\nabla^{\alpha}u_0(x)|\leq \frac{C(\alpha,u_0)}{|x|^{1+|\alpha|}}, ~~\forall ~~|\alpha|\ge 0.
\end{equation*}
Our first main result in this section is the following theorem.
\begin{theorem}{\rm (A-priori estimate for forward self similar solutions)}\label{th:aprioriEstimateForSelfSimilarSolution}\\
Let $u$, $u_0$ be as in the above. Then $U(\cdot):=u(\cdot,1)$, the solution profile at time $t=1$, belongs to $C^{\infty}(R^3)$ and
\begin{equation}
|\partial^{\alpha}\left(U(x)-e^{\Delta}u_0(x)\right)|\leq \frac{C(\alpha,u_0)}{(1+|x|)^{3+|\alpha|}}, ~~\forall~~|\alpha|\ge 0.
\end{equation}
\end{theorem}
\smallskip
\noindent
\textbf{Remarks:} Here and below, constants $C(u_0,\dots), T(u_0,\dots)\dots$ only depend on the magnitude of $u_0$ and its finitely many derivatives on the unit sphere. Similar estimates with more precise asymptotics have been proved in \cite{Brando} when the initial data is small in appropriate sense.\\

\noindent
\textbf{Proof:} Apply Lemma \ref{lm:aprioriEstimatesForLeraySolution} with $R=1$, we see (set $M:=\|u_0\|_{C(\partial B_1)}$)
\begin{equation}
 \sup_{0<t<T_1}\frac{1}{2}\int_{B_1(0)}|u(x,t)|^2dx+\int_0^{T_1}\int_{B_1(0)}|\nabla u(x,t)|^2dxdt\leq C(M),\quad T_1=T_1(M).
\end{equation}
For fixed $t_{\ast}<T_1$, with $t_{\ast}$ to be determined later, since $u(x,t)=\frac{1}{\sqrt{t}}u(\frac{x}{\sqrt{t}},1)=\frac{1}{\sqrt{t}}U(\frac{x}{\sqrt{t}})$, we have
\begin{eqnarray}\label{eq:localEnergyEstimate}
 C(M)&\ge&1/2\int_{B_1(0)}|u(x,t_{\ast})|^2dx+\int_{t_{\ast}/2}^{t_{\ast}}\int_{B_1(0)}|\nabla u(x,t)|^2dxdt\nonumber\\
&\ge&\frac{\sqrt{t_{\ast}}}{2}\int_{B_{\frac{1}{\sqrt{t_{\ast}}}}(0)}|u(x,1)|^2dx+\frac{\sqrt{t_{\ast}}}{8}\int_{B_{\frac{1}{\sqrt{t_{\ast}}}}(0)}|\nabla u(x,1)|^2dx.\\
&\ge&\frac{\sqrt{t_{\ast}}}{2}\int_{B_{\frac{1}{\sqrt{t_{\ast}}}}(0)}|U(x)|^2dx+\frac{\sqrt{t_{\ast}}}{8}\int_{B_{\frac{1}{\sqrt{t_{\ast}}}}(0)}|\nabla U(x)|^2dx.
\end{eqnarray}
On the other hand, for $\forall x_0$, $|x_0|=8$, since $u_0\in C^{\infty}(B_4(x_0))$, we can apply Theorem \ref{th:localRegularityForLeraySolution} and some simple boostraping arguments to show the following:\\
there exists $T_2=T_2(M)>0$ such that $\forall~\alpha$,
\begin{equation*}
\|\partial_t\partial_x^{\alpha}u\|_{L^{\infty}(B_{1/8}(x_0)\times[0,T_2])}\leq C(\alpha,u_0),
\end{equation*}
this is true for any $u\in\mathcal{N}(u_0)$.\\
Since $\forall~\lambda>0$, $\lambda u(\lambda x,\lambda^2 t)$ is also a Leray solution with initial data $u_0$, we obtain
\begin{equation*}
|\lambda^{|\alpha|+1}\partial^{\alpha}u(\lambda x_0,\lambda^2 t)-\partial^{\alpha}u_0(x_0)|\leq C(\alpha,u_0)t,
\end{equation*}
for any $\lambda>0$, $|\alpha|\ge 0$, $t\leq T_2(u_0)$.\\
Take $\lambda=\frac{1}{\sqrt{t}}$, we obtain $|(\frac{1}{\sqrt{t}})^{|\alpha|+1}\partial^{\alpha}u(\frac{x_0}{\sqrt{t}},1)-\partial^{\alpha}u(x_0)|\leq C(\alpha,u_0)t$.\\
Setting $y=\frac{x_0}{\sqrt{t}}$, and using the homogeneity of $\partial^{\alpha}u_0$, we get
\begin{equation}
 |\partial^{\alpha}U(y)-\partial^{\alpha}u_0(y)|\leq \frac{C(\alpha,u_0)}{|y|^{|\alpha|+3}},\quad\forall~~|y|>\frac{8}{\sqrt{T_2}}.
\end{equation}
Now choose $t_{\ast}$ sufficiently small, $t_{\ast}=t_{\ast}(M)$, we see from inequality (\ref{eq:localEnergyEstimate}):
$$\int_{B_{\frac{16}{\sqrt{T_2}}}}(|U(y)|^2+|\nabla U(y)|^2)dy\leq C(M).$$
Since $u(x,t)$ satisfies Navier-Stokes equations, it is easy to verify $U$ satisfies
\begin{eqnarray}
 \left.\begin{array}{rl}
        -\Delta U-\frac{x}{2}\cdot\nabla U-\frac{U}{2}+U\cdot\nabla U+\nabla P&=0\\
       \mbox{div~~}U&=0
       \end{array}\right\}\quad \mbox{in $R^3$.}
\end{eqnarray}
Thus elliptic estimates give
\begin{equation*}
\|U(\cdot)\|_{C^k(B_{\frac{9}{\sqrt{T_2}}})}\leq C(k,M).
\end{equation*}
These estimates, combined with the properties of heat equation, finish the proof.\\

For later use, let us study a nonhomogeneous Stokes system with singular forcing. Our result is the following Lemma.
\begin{lemma}{\rm (Decay for the linear singularly forced Stokes system)}\label{th:decayForSingularStokes}\\
Let $f\in C(R^3)$, suppose $v\in L^{\infty}_tL^{\gamma}_x(R^3\times(0,T))$ for any $T<\infty$, and some $\gamma>1$, suppose $v$ satisfies
\begin{eqnarray}
 \left.\begin{array}{rl}
        \partial_tv-\Delta v+\nabla p&=t^{-3/2}f(\frac{x}{\sqrt{t}})\\
       {\rm div~~}v&=0
       \end{array}\right\}\quad\mbox{in $R^3\times(0,\infty)$,}
\end{eqnarray}
for some distribution $p$, and $\lim_{t\to 0+}\|v(\cdot,t)\|_{L^{\gamma}(R^3)}=0$. Then\\
i) if $\tilde{v}$ also satisfies the above conditions, then $v=\tilde{v}$.\\
ii) if $f$ satisfies $M:=\sup_{x\in R^3}(1+|x|)^3|f(x)|<\infty$, then
\begin{equation}
v(\cdot,t)=\int_0^te^{\Delta(t-s)}P\frac{1}{s^{3/2}}f(\frac{\cdot}{\sqrt{s}})ds,
\end{equation}
where $P$ is the Helmholtz projection operator. Let $V(x)=v(x,1)$, then $\|V\|_{C^{1,\alpha}(B_R)}\leq C(\alpha,R)M$ for $\alpha\in (0,1)$ and
\begin{equation}
 \sup_{x\in R^3}\left((1+|x|)^2|V(x)|+(1+|x|)^3|\nabla V(x)|\right)\leq CM.
\end{equation}
iii) if $f$ satisfies $M:=\sup_{x\in R^3}(1+|x|)^4|f(x)|<\infty$, then
\begin{equation}
v(\cdot,t)=\int_0^te^{\Delta(t-s)}P\frac{1}{s^{3/2}}f(\frac{\cdot}{\sqrt{s}})ds.
\end{equation}
Let $V(x)=v(x,1)$, then $\|V\|_{C^{1,\alpha}(B_R)}\leq C(\alpha,R)M$ for $\alpha\in (0,1)$ and
\begin{equation}
\sup_{x\in R^3}\left((1+|x|)^{3}|V(x)|+(1+|x|)^{4}|\nabla V(x)|\right)\leq CM.
\end{equation}
\end{lemma}
\smallskip
\noindent
\textbf{Proof:} The uniqueness is easy. We only need to show that if $f=0$ and for some $\gamma_1,~\gamma_2>1$,
\begin{eqnarray*}
&& v\in L^{\infty}_t(L_x^{\gamma_1}+L_x^{\gamma_2})(R^3\times(0,T)) ~~{\rm for~~any~~} T>0 {\rm ~~and,~~}\\
&&\lim_{t\to 0+}\|v(t,\cdot)\|_{(L_x^{\gamma_1}+L_x^{\gamma_2})(R^3)}=0,
\end{eqnarray*}
then $v=0$. Set $\omega={\rm curl~~}v$, then $\partial_t\omega-\Delta\omega=0$ in $R^3\times(0,\infty)$. Since
\begin{equation*}
\lim_{t\to 0+}\|v(\cdot,t)\|_{(L^{\gamma_1}+L^{\gamma_2})(R^3)}=0,
\end{equation*}
we can extend $\omega$ to $R^3\times R$ by setting $\omega=0$ for $v<0$, and the extended function, which we still denote as $\omega$, satisfies $\partial_t\omega-\Delta \omega=0$ in $R^3\times R$. Here again there is no problem showing that the equation is satisfied across $t=0$ since $\omega$ decays to $0$ as $t\to 0+$. One can for example first mollify $\omega$ in $x$ and the mollified function is smooth in both $x$ and $t$. For the mollified function the claim is clear, then we can pass to the limit to show our claim. Since we have bounds for $\omega$ in some negative Sobolev space and $\omega=0$ for $t<0$, we conclude $\omega\equiv 0$. Thus $\Delta v=0$ in $R^3\times(0,\infty)$. Therefore $v\equiv 0$.\\
Let us now prove part ii) and part iii). By the uniqueness result, we only need to prove the claimed estimates. Denote the kernel of $Pe^{\Delta}$ by $k(x)$, then $k(\cdot)\in L^{1+\epsilon}(R^3)$ for any $\epsilon>0$. By Young's inequality it is easy to get
\begin{eqnarray*}
&&\|\int_0^te^{\Delta(t-s)}Ps^{-3/2}f(\frac{\cdot}{\sqrt{s}})ds\|_{L_x^{\frac{1+\epsilon}{1-\epsilon}}}\\
&&\leq C(\epsilon)\int_0^t\|(t-s)^{-3/2}k(\frac{\cdot}{\sqrt{t-s}})\|_{L^{1+\epsilon}_x}\|s^{-3/2}f(\frac{\cdot}{\sqrt{s}})\|_{L^{1+\epsilon}_x}ds\\
&&\leq C(\epsilon)Mt^{1-\frac{3\epsilon}{1+\epsilon}}.
\end{eqnarray*}
Thus,
\begin{equation*}
v(\cdot,t)=\int_0^te^{\Delta(t-s)}P\frac{1}{s^{3/2}}f(\frac{\cdot}{\sqrt{s}})ds.
\end{equation*}
Now let us prove the decay estimates of $V$. The proof is a direct consequence of the following inequality (which can be proved by simple calculations) with $\alpha,~\beta=3,~4$ and $R:=|x|>8$ :
\begin{equation}\label{eq:EstimateForAnIntegral}
 \int_0^1\int_{R^3}\frac{1}{(|x-y|+\sqrt{1-t})^{\alpha}}\frac{1}{(|y|+\sqrt{t})^{\beta}}dydt\leq \left\{\begin{array}{ll}
                                                                                              R^{-3}\log R &\quad {\rm if~~} \alpha=\beta=3,\\
                                                                                             R^{-\alpha-\beta+4}&\quad {\rm otherwise}.
                                                                                             \end{array}\right.
\end{equation}
For part i), we have
\begin{eqnarray*}
 &&|V(x)|\leq\int_0^1\int_{R^3}\frac{1}{(|x-y|+\sqrt{1-t})^{3}}\frac{1}{(|y|+\sqrt{t})^3}dydt\leq |x|^{-3}\log |x|,\\
&&|\nabla V(x)|\leq \int_0^1\int_{R^3}\frac{1}{(|x-y|+\sqrt{1-t})^{4}}\frac{1}{(|y|+\sqrt{t})^3}dydt\leq |x|^{-3},
\end{eqnarray*}
for $|x|>8$. For part ii), we have
\begin{eqnarray*}
 &&|V(x)|\leq\int_0^1\int_{R^3}\frac{1}{(|x-y|+\sqrt{1-t})^{3}}\frac{1}{(|y|+\sqrt{t})^4}dydt\leq |x|^{-3},\\
&&|\nabla V(x)|\leq \int_0^1\int_{R^3}\frac{1}{(|x-y|+\sqrt{1-t})^{4}}\frac{1}{(|y|+\sqrt{t})^4}dydt\leq |x|^{-4},
\end{eqnarray*}
for $|x|>8$. Thus the decay estimates are proved.
Since $V$ also satisfies an elliptic equation:
\begin{eqnarray}
 \left.\begin{array}{rl}
        -\Delta V-\frac{x}{2}\cdot\nabla V-\frac{V}{2}+\nabla P&=f\\
            \mbox{div~~}V&=0
       \end{array}\right\}\quad\mbox{in $R^3$,}
\end{eqnarray}
the estimates in $B_C(0)$ is simple.\\
\end{section}

\begin{section}{Existence of forward self similar solution for large initial data}
In this section, we prove the following theorem.
\begin{theorem}\label{th:existenceForSmooth}
Let $u_0\in C^{\infty}(R^3\setminus\{0\})$ satisfy $\lambda u_0(\lambda x)=u_0(x)$ for all $\lambda>0$, ${\rm div~~}u_0=0$. Then there exists $u\in C^{\infty}(R^3\times(0,\infty))$, with $\lambda u(\lambda x,\lambda^2 t)=u(x,t)$ for all $\lambda>0$, and $u\in \mathcal{N}(u_0)$, that is, $u$ satisfies
\begin{eqnarray}
 \left.\begin{array}{rl}
        \partial_tu-\Delta u+u\cdot\nabla u+\nabla p&=0\\
        {\rm div~~}u&=0
       \end{array}\right\}\quad\mbox{in $R^3\times(0,\infty)$ for some $p$. }
\end{eqnarray}
Moreover, let $U(x)=u(x,1)$, then
\begin{equation*}
|\partial^{\alpha}\left(U(x)-e^{\Delta}u_0(x)\right)|\leq \frac{C(\alpha,u_0)}{(1+|x|)^{3+|\alpha|}}.
\end{equation*}
\end{theorem}
\smallskip
\noindent
\textbf{Proof:} By Theorem \ref{th:aprioriEstimateForSelfSimilarSolution}, it suffices to show there exists $u\in \mathcal{N}(u_0)$ with the scaling
\begin{equation}
\lambda u(\lambda x,\lambda^2t)=u(x,t) {\rm ~~~for~~ all~} \lambda>0.
\end{equation}
Denote
\begin{equation}
X=\{V\in C^1(R^3): \mbox{div~~}V=0, ~~\sup_{x\in R^3}\left((1+|x|)^2|V(x)|+(1+|x|)^3|\nabla V(x)|\right)<\infty\}.
\end{equation}
For any $V\in X$, we define a natural norm
\begin{equation}
\|V\|_X=\sup_{x\in R^3}\left((1+|x|)^2|V(x)|+(1+|x|)^3|\nabla V(x)|\right).
\end{equation}
Set $U_0=e^{\Delta}u_0$. Introduce a parameter $\mu\in[0,1]$, set $U_{0\mu}=\mu U_0$. We will follow Leray's method to prove the existence of $u\in \mathcal{N}(u_0)$ with $\lambda u(\lambda x,\lambda^2 t)=u(x,t)$ for all $\lambda>0$. Due to the scaling invariance of $u(x,t)$, we are essentially seeking the profile function $U(x)=u(x,1)$, where $U(x)$ satisfies
\begin{eqnarray}
 \left.\begin{array}{rl}
        -\Delta U+U\cdot\nabla U-\frac{U}{2}-\frac{x}{2}\cdot\nabla U+\nabla P&=0\\
         \mbox{div~~}U&=0
       \end{array}\right\}\quad\mbox{in $R^3$,}
\end{eqnarray}
and the correct asymptotics at spatial infinity. We will solve $U$ in the following form
\begin{equation}
U=U_{0\mu}+V, {\rm ~~where~~} V\in X.
\end{equation}
It is clear $u(x,t)=\frac{1}{\sqrt{t}}U(\frac{x}{\sqrt{t}})\in\mathcal{N}(\mu u_0)$ if and only if $U(x)$ satisfies the above elliptic system and $U(x)=U_{0\mu}+V$ for some $V\in X$, by Theorem \ref{th:aprioriEstimateForSelfSimilarSolution}. Thus we have reduced the problem to to finding $V\in X$, with
\begin{eqnarray}
 \left.\begin{array}{rl}
-\Delta V+V\cdot\nabla V+U_{0\mu}\cdot\nabla V+V\cdot\nabla U_{0\mu}-\frac{V}{2}-\frac{x}{2}\cdot\nabla V+\nabla P&=-U_{0\mu}\cdot\nabla U_{0\mu}\\
         \mbox{div~~}V&=0
       \end{array}\right\},
\end{eqnarray}
in $R^3$. We rewrite the above as:
\begin{equation}
-\Delta V-\frac{V}{2}-\frac{x}{2}\cdot\nabla V+\nabla P=-V\cdot\nabla V-U_{0\mu}\cdot\nabla V-V\cdot\nabla U_{0\mu}-U_{0\mu}\cdot\nabla U_{0\mu}.
\end{equation}
Since $V\in X$, $V$ satisfies the above equation if and only if $v(x,t):=\frac{1}{\sqrt{t}}V(\frac{x}{\sqrt{t}})$ satisfies
 \begin{eqnarray}\label{eq:singularStokesSolution}
 \left.\begin{array}{rl}
        \partial_tv-\Delta v+\nabla p&=t^{-3/2}F(\frac{x}{\sqrt{t}})\\
       \mbox{div~~}v&=0\\
      v(\cdot,0)&=0
       \end{array}\right\}
\end{eqnarray}
where
\begin{equation}
F=-V\cdot\nabla V-U_{0\mu}\cdot\nabla V-V\cdot\nabla U_{0\mu}-U_{0\mu}\cdot\nabla U_{0\mu}
\end{equation}
has the decay properties in Lemma \ref{th:decayForSingularStokes}. Thus for such $F$, equation (\ref{eq:singularStokesSolution}) is uniquely solvable, we denote the solution profile at time $1$ as $\mathcal{G}(F)\in X$. This enables us to consider the following equivalent formulation,
\begin{equation}
\mbox{find $V\in X$ with~~} V=\mathcal{G}(-V\cdot\nabla V-U_{0\mu}\cdot\nabla V-V\cdot\nabla U_{0\mu}-U_{0\mu}\cdot\nabla U_{0\mu})\mbox{.}
\end{equation}
Let $K:X\times[0,1]\rightarrow X$ be defined as:
\begin{equation}
\forall~~V\in X,~~\mu\in[0,1],~~K(V,\mu):=\mathcal{G}(U_{0\mu}\nabla U_{0\mu})+\mathcal{G}(U_{0\mu}\nabla V+V\nabla U_{0\mu}+V\nabla V).
\end{equation}
The first term on the right hand side is one dimensional. The second term by estimates in Lemma \ref{th:decayForSingularStokes} is compact. The compactness is due to the local $C^{1,\alpha}$ estimates and the fast decay at inifinity. Thus we conclude $K\in C^1(X\times[0,1])$ is compact. Therefore we are reduced to solve the following abstract problem:\\
\begin{center}
find $V\in X$, such that $V+K(V,\mu)=0$,  where $\mu\in[0,1]$.                                                         \end{center}
\bigskip
At this stage, we are in a position to apply Leray's method, see for example \cite{Maw}. We need the following conditions to be verified:\\
1. Solvability for $\mu$ small. This is already done, for example in \cite{CaPl,GiMi}, note that it also follows from a simple implicit function theorem in our formulation. In the language of Leray Schauder degree theory, we can verify $d(I+K(\cdot,\mu),B_M(0),0)=1$ for $\mu$ small and some fixed $M>0$.\\
2. A priori estimate for solutions. This is done, in Theorem \ref{th:aprioriEstimateForSelfSimilarSolution}.\\
3. Compactness and continuity of $K$. This follows from the estimates of $\mathcal{G}$.\\
Thus we can apply Leray's method, and conclude that for each $\mu\in [0,1]$, there exists a solution $V\in X$ to $V+K(V,\mu)=0$. Take $\mu=1$, the theorem is proved. \\

With the existence theorem for smooth (away from $0$) $-1$ homogeneous initial data, we can obtain existence results for not so smooth initial data. We illustrate the method with H\"{o}lder continuous (away from $0$) initial data, although more general initial data can be considered.
\begin{theorem}
Let $u_0\in C_{\rm loc}^{\alpha}(R^3\setminus \{0\})$ with $\alpha\in(0,1)$, $\lambda u_0(\lambda x)=u_0(x)$ for all $\lambda>0$, and ${\rm ~div~~}u_0=0$ in $R^3$. Denote $M=\|u_0\|_{C^{\alpha}(\partial B_1)}$. Then there exists $u\in\mathcal{N}(u_0)$, and $u$ satisfies $u(x,t)=\lambda u(\lambda x,\lambda^2t)$ for all $\lambda>0$. Moreover, let $U(x)=u(x,1)$. Then $U\in C^{\infty}(R^3)$ with
\begin{equation}
 |U(x)-e^{\Delta}u_0(x)|\leq \frac{C(M)}{(1+|x|)^{1+\alpha}}.
\end{equation}
\end{theorem}
\smallskip
\noindent
\textbf{Proof:} Let us choose $u_0^{\epsilon}\in C^{\infty}(R^3\setminus \{0\})$ with $\lambda u_0^{\epsilon}(\lambda x)=u_0^{\epsilon}(x)$ for all $\lambda>0$, ${\rm ~div~~}u_0^{\epsilon}=0$ in $R^3$, $\|u_0^{\epsilon}\|_{C^{\alpha}(\partial B_1(0))}\leq CM$, and $\|u_0^{\epsilon}-u_0\|_{C(\partial B_1)}\to 0$ as $\epsilon \to 0+$. We can construct such $u_0^{\epsilon}$ by first mollifying $u_0$ on the unit sphere and then using the scaling invariance and applying Helmholtz projection operator to form $u_0^{\epsilon}$. We only note that the scaling invariance is preserved by the Helmholtz projection. By Theorem \ref{th:existenceForSmooth}, we can find $u^{\epsilon}\in \mathcal{N}(u_0^{\epsilon})$ with $\lambda u^{\epsilon}(\lambda x,\lambda^2t)=u^{\epsilon}(x,t)$, for all $\lambda>0$. Let $U^{\epsilon}(x)=u^{\epsilon}(x,1)$, then $u^{\epsilon}(x,t)=\frac{1}{\sqrt{t}}U^{\epsilon}(\frac{x}{\sqrt{t}})$. For any $x_0\in R^3$ with $|x_0|=8$, since $u_0^{\epsilon}\in C^{\alpha}(B_4(x_0))$ with $\|u_0^{\epsilon}\|_{C^{\alpha}(B_4(x_0))}\leq C(M)$, by Theorem \ref{th:localHolderLeraySolution}, there exists $T(M)>0$, such that $u^{\epsilon}\in C_{{\rm par}}^{\alpha}(B_{1/2}\times [0,T(M)])$ and $\|u^{\epsilon}\|_{C_{{\rm par}}^{\alpha}(B_{1/2}\times[0,T(M)])}\leq C(M)$. Thus,
\begin{equation}
 |\frac{1}{\sqrt{t}}U^{\epsilon}(\frac{x_0}{\sqrt{t}})-u_0^{\epsilon}(x_0)|\leq C(M)t^{\alpha/2}, \quad{\rm for ~~} t<T(M).
\end{equation}
By the homogeneity of $u_0^{\epsilon}$, we get
\begin{equation}
 |U^{\epsilon}(\frac{x_0}{\sqrt{t}})-u_0^{\epsilon}(\frac{x_0}{\sqrt{t}})|\leq C(M)t^{1/2+\alpha/2},\quad{\rm for~~} t<T(M).
\end{equation}
Notice that $|x_0|=8$ is arbitrary, we get
\begin{equation}
 |U^{\epsilon}(x)-u_0^{\epsilon}(x)|\leq \frac{C(M)}{|x|^{1+\alpha}}\quad{\rm for~~} |x|>C_1(M).
\end{equation}
 Moreover, by following the same arguments in the proof of Theorem \ref{th:aprioriEstimateForSelfSimilarSolution}, we can obtain
\begin{equation}
 \|U^{\epsilon}\|_{C^k(B_R(0))}\leq C(k,M,R) \quad {\rm for~~} \forall R>0.
\end{equation}
By combining the above estimates and using elementary properties of heat equation, we get
\begin{equation}
 |U^{\epsilon}(x)-e^{\Delta}u^{\epsilon}(x)|\leq \frac{C(M)}{(1+|x|)^{1+\alpha}}, \quad {\rm for ~~}x\in R^3.
\end{equation}
Note also that since $u^{\epsilon}$ satisfies the Navier-Stokes equations for $t>0$, $U^{\epsilon}$ satisfies
\begin{eqnarray*}
 \left.\begin{array}{rl}
        -\Delta U^{\epsilon}+ U^{\epsilon}\cdot\nabla  U^{\epsilon}-\frac{x}{2}\cdot\nabla  U^{\epsilon}-\frac{ U^{\epsilon}}{2}+\nabla P^{\epsilon}&=0\\
    {\rm div~~} U^{\epsilon}&=0
       \end{array}\right\}{\rm in~~} R^3.
\end{eqnarray*}
By the estimates on $U^{\epsilon}$, we can pass to a subsequence $\epsilon_i\to 0+$, such that $U^{\epsilon_i}\to U$ in $C^2(B_R(0))$ for all $R>0$. Thus $U$ satisfies
\begin{eqnarray}
 \left.\begin{array}{rl}
        -\Delta U+ U\cdot\nabla U-\frac{x}{2}\cdot\nabla  U-\frac{U}{2}+\nabla P&=0\\
    {\rm div~~} U&=0
       \end{array}\right\}{\rm in~~} R^3,
\end{eqnarray}
and
\begin{equation}
 |U(x)-e^{\Delta}u_0(x)|\leq \frac{C(M)}{(1+|x|)^{1+|\alpha|}}~~{\rm for~~all~~}x\in R^3.
\end{equation}
Setting $u(x,t)=\frac{1}{\sqrt{t}}U(\frac{x}{\sqrt{t}})$, we can easily verify that $u$ satisfies all the conditions in our theorem.\\

\end{section}

\bigskip
\noindent
{\bf Acknowledgements.}
We thank Gregory Seregin for valuable comments.\\
\noindent
This work was supported in part by grant
 DMS 1101428 from the National Science Foundation.

\end{section}

\end{document}